\documentclass[12pt,a4paper]{article}
\usepackage[utf8]{inputenc}
\usepackage{amssymb, amsthm, amsmath}
\usepackage[english]{babel}
\usepackage[T1]{fontenc}
\usepackage{graphicx}
\usepackage{mathtools}
\usepackage[font=small,labelfont={normal, bf}]{caption}
\usepackage{subcaption}
\usepackage{xcolor}

\usepackage[cm]{fullpage}

\usepackage{autonum}

\makeatletter
\newcommand{\@giventhatstar}[2]{#1\;\middle|\;#2)}
\newcommand{\@giventhatnostar}[3][]{#1#2\;#1|\;#3#1}
\newcommand{\giventhat}{\@ifstar\@giventhatstar\@giventhatnostar}
\makeatother

\newtheorem{theorem}{Theorem}

\theoremstyle{definition}
\newtheorem{definition}{Definition}

\newtheorem{lemma}{Lemma}
\DeclareMathOperator{\sgn}{sgn}

\title{Scaling limits of Lévy walks with random velocities}
\author{
Hubert Woszczek\thanks{Faculty of Pure and Applied Mathematics, Wroclaw University of Science and Technology, Wyb. Wyspia\'nskiego 27, 50-370 Wroc{\l}aw, Poland}$^{*,}$\thanks{\underline{Corresponding author:} \texttt{hubert.woszczek@pwr.edu.pl}}, \and Marek A. Teuerle$^*$, \and Agnieszka Wy\l{}oma\'nska$^*$
}
\date{}

\begin{document}

\maketitle
\begin{abstract}
    This paper investigates Lévy walks with random velocities, extending classical models beyond constant speed assumptions. We derive scaling limits, demonstrating that diffusion depends on interplay between heavy-tailed duration and velocity distributions. Three distinct scaling regimes are identified, including a critical case with logarithmic corrections, offering a precise framework for modeling anomalous transport in heterogeneous systems.
\end{abstract}



\section{Introduction}
The mathematical modeling of transport phenomena has been significantly enriched by the introduction of the Lévy walk (LW) framework by Shlesinger, Klafter, and Wong in 1982 \cite{Shlesinger1982}. Unlike the standard Brownian motion, the LW model describes a class of movements that exhibit anomalous diffusion properties through a specific mechanism: a particle moves at a constant velocity during segments whose durations are drawn from a heavy-tailed distribution. This construction addresses the physical inconsistencies of Lévy flights, which assume instantaneous displacements and thus infinite propagation speeds. By enforcing a strict spatiotespamporal coupling, the LW ensures that the particle's displacement is explicitly bounded by its travel time.

The practical relevance of the LW framework has been demonstrated in various fields, such as the dynamics of cold atoms in optical lattices \cite{Sagi2012}, swarming bacteria migration \cite{Ariel2015}, and the optimized foraging strategies observed in diverse biological organisms \cite{Viswanathan1999}. A comprehensive systematic categorization of these models is provided in the review by Zaburdaev, Denisov, and Klafter \cite{Zaburdaev2015}. Their work highlights a wide spectrum of LW variants and illustrates how different forms of spatiotemporal coupling result in distinct scaling regimes that deviate from the classic Gaussian framework. Building upon these foundations, recent research has expanded into more complex scenarios, including multidimensional environments and heterogeneous media. A significant contribution to this evolution is the work of Magdziarz and Teuerle \cite{Magdziarz2015}, who provided a detailed analysis of the asymptotic properties and numerical simulations of multidimensional Lévy walks. Their findings showed that the scaling limits of such processes depend critically on the power-law exponents of the underlying distributions, often converging to subordinated stable motions or Brownian motion depending on the specific coupling regimes.

Furthermore, the theoretical framework has been deepened by limit theorems that provide the governing equations, involving the so-called fractional material derivatives  \cite{sokolov2003towards,Plociniczak2024}, for these complex processes \cite{Magdziarz2015}. This article extends these discussions by focusing on the dynamics of Lévy walks with random velocities. In many disordered systems, the assumption of constant speed is insufficient to capture the inherent variability of the environment \cite{Zaburdaev2015}. By treating the velocity of each step as a stochastic variable, we investigate the interplay between heavy-tailed waiting times and heavy-tailed velocity distributions. This analysis identifies critical thresholds at which velocity fluctuations dominate transport dynamics, providing a more robust framework for predicting the long-term behavior of particles in fluctuating environments.

\section{Preliminaries}
Let us begin with the notation convention used in this article. For two functions $f,g$, the notion $f\sim g$ means that $\lim_{x\to x_0}f(x)/g(x)=1$. $f\lesssim g$ means that there exists a function $h$ such that $f\le h$ and $h\sim g$. We will also use the big O notation $f(x) = \mathcal{O}(g(x))$ as $x\to x_0$, which means that there exists a positive constant $M$ such that $|f(x)|\le Mg(x)$ on the interval $[x_0-y, x_0+y]$ for some constant $y>0$. If $x_0$ is 'equal' to $\infty$, we mean the interval $[x_0, \infty]$. Convergence in finite-dimensional distributions is denoted by $\overunderset{d}{n \to \infty}{\to }$, and convergence in Skorokhod topology is denoted by $\overunderset{\mathbb{J}_1}{n \to \infty}{\to }$ (for more details, see \cite{billingsley1999convergence}).

\noindent Next, we present the basic notation and facts about general continuous time random walks (CTRWs). Let $\left\{T_i\right\}_{i \in \mathbb{N}_+}$ be the sequence of positive independent identically distributed (IID) random variables representing the waiting times between consecutive jumps of the walker. The counting process
\begin{equation}
    N\left(t\right) = \max\left\{k \geq 0: \sum_{i=1}^k T_i \leq t\right\}
\end{equation}
counts the number of jumps of the walker up to time t. It is so-called a wait-first scenario. Next, let $\left\{\mathbf{J}_i\right\}_{i \in \mathbb{N}_+}$ be the sequence of IID random vectors that represent the consecutive jumps of the walker. Then, the CTRW
\begin{equation}
    R\left(t\right) = \sum_{i=1}^{N\left(t\right)} \mathbf{J}_i
\end{equation}
describes the position of the walker at time $t$.

Let us now assume that waiting times are heavy-tailed distributed random variables with index $\alpha \in \left(0, 1\right)$, i.e.
\begin{equation}\label{heavytail}
    \mathbb{P}\left(T_i>t\right) \sim ct^{-\alpha}
\end{equation}
for some positive constant $c$ (we can assume $c=1$ for simplicity). Then, the sequence $\{T_i\}_{i\in\mathbb{N}_+}$ of waiting times belongs to the domain of attraction of one-sided $\alpha$-stable distribution, namely
\begin{equation}\label{scallimsalpha}
    n^{-1/\alpha}\sum_{i=1}^{\left[nt\right]} T_i \overunderset{d}{n \to \infty}{\to } S_{\alpha}\left(t\right),
\end{equation}
for fixed $t$, where $\{S_{\alpha}\left(t\right)\}_{t\ge0}$ ($\alpha \in \left(0, 1\right)$) is a $\alpha$-stable subordinator with the Laplace transform
\begin{equation}
    \mathbb{E}\left[\exp\left\{-sS_{\alpha}\left(t\right)\right\}\right] = \exp\left\{-ts^{\alpha}\right\}.
\end{equation}
Its Fourier transform is given by
\begin{equation}
    \mathbb{E}\left[\exp\left\{ikS_{\alpha}\left(t\right)\right\}\right] = \exp\left\{-t\left|k\right|^{\alpha}\cos\left(\pi\alpha/2\right)\left(1- i\tan\left(\pi\alpha/2\right)\sgn\left(k\right)\right)\right\}.
\end{equation}
This process is uniquely determined by its L\'evy triplet
\begin{equation}
    \left[\int_{x\neq0}\frac{x}{1+x^2}\nu_{\alpha}\left(dx\right), 0, \nu_{\alpha}\left(dx\right)=\frac{\alpha x^{-\alpha-1}}{\Gamma\left(1-\alpha\right)}\mathbf{1}_{\left\{x>0\right\}}dx\right].
\end{equation}
We will also recall the definition of the so-called inverse $\alpha$-stable subordinator $\{S_{\alpha}^{-1}\left(t\right)\}_{t\ge0}$, defined by
\begin{equation}
    S_{\alpha}^{-1}\left(t\right) = \inf\left\{\tau \geq 0: S_{\alpha}\left(\tau\right)>t\right\}.
\end{equation}
The process $\{S_{\alpha}^{-1}\left(t\right)\}_{t\ge0}$ is the scaling limit of the counting process $\{N\left(t\right)\}_{t\ge0}$ with $\{T_i\}_{i\in\mathbb{N}_0}$ satisfying the condition \eqref{heavytail}. Additionally, we assume that jumps of the walker are equal to
\begin{equation}\label{jumpdef}
    \mathbf{J}_i = V_i \mathbf{I}_i T_i,
\end{equation}
where $\left\{V_i\right\}_{i \in \mathbb{N}}$ is a sequence of IID random variables corresponding to the walker velocities, independent of $\mathbf{I}_i$, $T_i$ for all $i \in \mathbb{N}$ and fulfills condition \eqref{heavytail} with the stability index $\beta$ for some $\beta \in \left(0, 1\right)$. $\left\{I_i\right\}_{i \in \mathbb{N}}$ is a sequence of non-degenerate unit vectors from $\mathbb{R}^d$ corresponding to the walker directions, independent of $V_i$ and $T_i$ for all $i \in \mathbb{N}$.
\begin{definition}\label{lwrv}
    Let $\left\{\left(V_i, T_i, \mathbf{J}_i\right)\right\}_{i\geq1}$ be the sequence of IID random vectors with components defined in Eq. \eqref{jumpdef}, then
    \begin{equation}
        U_V\left(t\right) = \sum_{i=1}^{N\left(t\right)} \mathbf{J}_i, \quad O_V\left(t\right) = \sum_{i=1}^{N\left(t\right)+1} \mathbf{J}_i, \quad W_V\left(t\right) = \sum_{i=1}^{N\left(t\right)} \mathbf{J}_i + \left(t - \sum_{i=1}^{N\left(t\right)}T_i\right)\mathbf{I}_{N\left(t\right)+1}V_{N\left(t\right)+1}
    \end{equation}
    are called wait-first L\'evy walk with random velocity, jump-first L\'evy walk with random velocity, and (continuous) L\'evy walk with random velocity, respectively. A similar model was introduced in the physics literature in \cite{Denisov2012Levy}.
\end{definition}
\section{Main results}
In this section, we will prove the scaling limits of L\'evy walks with random velocities. 
\noindent In \cite{Magdziarz2015} it was shown that
\begin{equation}\label{auxilaryjumpslimit}
    n^{-1/\alpha}\sum_{i=1}^{\left[nt\right]} \mathbf{I}_iT_i \overunderset{\mathbb{J}_1}{n \to \infty}{\to } \mathbf{L}_{\alpha}\left(t\right),
\end{equation}
where $\mathbf{L}_\alpha\left(t \right)$ is described by L\'evy triplet
\begin{equation}
    \left[\int_{\mathbf{x}\neq\mathbf{0}}\frac{\mathbf{x}}{1+\left\|\mathbf{x}\right\|^2}\nu_{\mathbf{L}_{\alpha}}\left(d\mathbf{x}\right), \mathbf{0}, \nu_{\mathbf{L}_{\alpha}}\left(d\mathbf{x}\right)=\int_{\mathbf{D}}\frac{\alpha r^{-\alpha-1}}{\Gamma\left(1-\alpha\right)}\mathbf{1}_{\left\{r>0\right\}}\mathbf{\Lambda}\left(d\mathbf{u}\right)\right].
\end{equation}
In the above, $\mathbf{x}\in \mathbb{R}^d\setminus\left\{\mathbf{0}\right\}$, $r=\left\|\mathbf{x}\right\|$, $\mathbf{u} = \mathbf{x}/\left\|\mathbf{x}\right\|\in \mathbf{D}\subset\mathbb{S}^{d-1}$ and $\mathbf{\Lambda}\left(d\mathbf{u}\right) = \mathbb{P}\left(\mathbf{I}_1\in d\mathbf{u}\right)$ is called the spectral measure. Next, we derive two auxiliary lemmas related to the scaling limit of our model.

\begin{lemma}\label{auxlemma1}
 Let $\left\{V_i\right\}_{i\geq 1}$, $\left\{T_i\right\}_{i\geq 1}$ be defined as in Eq. \eqref{jumpdef}. Then, we have
    \begin{equation}
    n^{-1/\alpha^*}\sum_{i=1}^{\left[nt\right]} V_iT_i \overunderset{\mathbb{J}_1}{n \to \infty}{\to } S_{\alpha^*}\left(t\right),
\end{equation}
where $S_{\alpha^*}\left(t \right)$ is described by L\'evy triplet
\begin{equation}
    \left[\int_{x\neq0}\frac{x}{1+x^2}\nu_{\alpha^*}\left(dx\right), 0, \nu_{\alpha^*}\left(dx\right)=\frac{\alpha x^{-\alpha^*-1}}{\Gamma\left(1-\alpha^*\right)}\mathbf{1}_{\left\{x>0\right\}}dx\right],
\end{equation}
where $\alpha^* = \min\left\{\alpha, \beta\right\}$.
\begin{proof}
    The distributions of $V_i$ and $T_i$ satisfy the condition \eqref{heavytail}, so their tails vary regularly with the indices $\beta$ and $\alpha$, respectively. From Theorem A of \cite{Kasahara2018} we know that the product of distributions with regularly varying tails is regularly varying with the index $\alpha^* = \min\left\{\alpha, \beta\right\}$. Thus, $V_iT_i$ belongs to the domain of attraction of one-sided L\'evy stable distribution. Applying \eqref{scallimsalpha} to our sequence $\left\{V_iT_i\right\}_{i\geq 1}$, we obtain the thesis.
\end{proof}
\end{lemma}
\begin{lemma}
    Let $\left\{\mathbf{J}_i\right\}_{i\geq 1}$ be defined as in Eq. \eqref{jumpdef}. Then, we have the following
    \begin{equation}
    n^{-1/\alpha^*}\sum_{i=1}^{\left[nt\right]} \mathbf{J}_i \overunderset{\mathbb{J}_1}{n \to \infty}{\to } \mathbf{L}_{\alpha^*}\left(t\right),
\end{equation}
where $\mathbf{L}_\alpha\left(t \right)$ is described by L\'evy triplet
\begin{equation}
    \left[\int_{\mathbf{x}\neq\mathbf{0}}\frac{\mathbf{x}}{1+\left\|\mathbf{x}\right\|^2}\nu_{L_{\alpha^*}}\left(d\mathbf{x}\right), \mathbf{0}, \nu_{L_{\alpha^*}}\left(d\mathbf{x}\right)=\int_{\mathbf{D}}\frac{\alpha^* r^{-\alpha^*-1}}{\Gamma\left(1-\alpha^*\right)}\mathbf{1}_{\left\{r>0\right\}}\mathbf{\Lambda}\left(d\mathbf{u}\right)\right],
\end{equation}
and $\alpha^* = \min\left\{\alpha, \beta\right\}$.
\begin{proof}
    Applying Lemma \ref{auxlemma1} and Eq. \eqref{auxilaryjumpslimit} to jumps, we obtain the thesis.
\end{proof}
\noindent The next result shows the behavior of the tails of the product of random variables $T$ and $V$, which represent waiting times and velocities respectively in the critical case $\alpha=\beta$.
\begin{lemma}
    Let $V, T$ be independent, positive heavy-tailed random variables with the same index of stability $\alpha$. Then
\begin{equation}
    \mathbb{P}(VT > z) \sim \alpha z^{-\alpha}\log z, \quad \text{as } z \to \infty.
\end{equation}
\begin{proof}
By the definition of asymptotic equivalence, for any arbitrarily small $\epsilon > 0$, there exists a sufficiently large $M > 0$ such that for all $x \ge M$, we have
\begin{equation}\label{eq:epsilon_bound}
    (1-\epsilon)x^{-\alpha} \le \bar{F}_V(x) \le (1+\epsilon)x^{-\alpha}.
\end{equation}
Using the law of total probability and the conditioning on $T = x$, the tail of the product is given by
\begin{equation}\label{eqint}
    H(t) = \mathbb{P}(T V > t) = \int_0^\infty \mathbb{P}(V > t/x) dF_T(x) = \int_0^\infty \bar{F}_V(t/x) dF_T(x).
\end{equation}
Assume $t$ is large enough such that $t > M^2$ (which implies $M < t/M$). We split the integral in Eq. \eqref{eqint} into three parts
\begin{equation}
    H(t) = \int_0^M \bar{F}_V(t/x) dF_T(x) + \int_M^{t/M} \bar{F}_V(t/x) dF_T(x) + \int_{t/M}^\infty \bar{F}_V(t/x) dF_T(x) = I_1(t) + I_2(t) + I_3(t).
\end{equation}
First, let us analyze $I_1(t)$. Since $x \in (0, M]$, thus we have $t/x \ge t/M \ge M$. Applying the upper bound from \eqref{eq:epsilon_bound}, we have
\begin{equation}
    I_1(t) \le \int_0^M (1+\epsilon)(t/x)^{-\alpha} \, dF_T(x) = (1+\epsilon)t^{-\alpha} \int_0^M x^\alpha \, dF_T(x).
\end{equation}
Because $\int_0^M x^\alpha \, dF_T(x)\le M^{\alpha}\int_0^M dF_T(x)\leq M^{\alpha}$, where $M^{\alpha}$ is a finite constant independent of $t$, thus $I_1(t) = \mathcal{O}(t^{-\alpha})$. For $I_3(t)$, we use the trivial bound that probability is not greater than $1$ and obtain
\begin{equation}
    I_3(t) = \int_{t/M}^\infty \mathbb{P}(V > t/x) \, dF_T(x) \le \int_{t/M}^\infty 1 \, dF_T(x) = \bar{F}_T(t/M).
\end{equation}
Since $t/M \ge M$, thus $\bar{F}_T(t/M) \sim (t/M)^{-\alpha} = M^\alpha t^{-\alpha}$. Therefore, $I_3(t) = \mathcal{O}(t^{-\alpha})$. \\
Since both $I_1(t)$ and $I_3(t)$ lack a logarithmic factor, they are negligible compared to $t^{-\alpha} \log(t)$ as $t \to \infty$. Now, we analyze $I_2(t)$. The integration limits ensure that both $x \ge M$ and $t/x \ge M$. We can therefore apply the bounds from \eqref{eq:epsilon_bound} to the integrand $\bar{F}_V(t/x)$
\begin{equation}
    (1-\epsilon) t^{-\alpha} \int_M^{t/M} x^\alpha \, dF_T(x) \le I_2(t) \le (1+\epsilon) t^{-\alpha} \int_M^{t/M} x^\alpha \, dF_T(x).
\end{equation}
Let $J(y) = \int_M^y x^\alpha \, dF_T(x)$ where $y = t/M$. Using integration by parts and noting that $dF_T(x) = -d\bar{F}_T(x)$, we have
\begin{equation}
\begin{aligned}
    J(y) &= -\int_M^y x^\alpha \, d\bar{F}_T(x) = -\left( \left[ x^\alpha \bar{F}_T(x) \right]_M^y - \int_M^y \bar{F}_T(x) \, d(x^\alpha) \right) \\
         &= M^\alpha \bar{F}_T(M) - y^\alpha \bar{F}_T(y) + \int_M^y \bar{F}_T(x) \alpha x^{\alpha-1} \, dx.
\end{aligned}
\end{equation}
As $y \to \infty$, $M^\alpha \bar{F}_T(M)$ is a constant, $y^\alpha \bar{F}_T(y) \to 1$ (since $\bar{F}_T(y) \sim y^{-\alpha}$) and the integral term dominates. Because $\bar{F}(x) \sim x^{-\alpha}$, the integrand is asymptotically $\alpha x^{-1}$. Thus, we have
    \begin{equation}
        \int_M^y \alpha x^{-1} \, dx = \alpha (\log(y) - \log(M)) \sim \alpha \log(y).
    \end{equation}
Therefore, $J(y) \sim \alpha \log(y)$. Substituting $y = t/M$ gives $J(t/M) \sim \alpha \log(t)$. Substituting the asymptotic behavior of $J(t/M)$ back into our bounds for $I_2(t)$ yields
\begin{equation}
    (1-\epsilon) \alpha t^{-\alpha} \log(t) \lesssim I_2(t) \lesssim (1+\epsilon) \alpha t^{-\alpha} \log(t).
\end{equation}
Since $I_1(t)$ and $I_3(t)$ are $\mathcal{O}(t^{-\alpha})$, the total probability $H(t)$ is asymptotically equivalent to $I_2(t)$. Because $\epsilon$ was arbitrarily chosen, we conclude that
$   \lim_{t \to \infty} {\mathbb{P}(T V > t)}/{\alpha t^{-\alpha} \log(t)} = 1.$
This completes the proof.
\end{proof}
\end{lemma}
\end{lemma}
\begin{lemma}\label{vectorconv}
    The joint vector of partial sums of jumps and waiting times has the following convergence in the Skorokhod topology:
    \begin{enumerate}
        \item
        for $\alpha<\beta$, we have
        \begin{equation}\label{convcase1}
            \left(n^{-1/\alpha}\sum_{i=1}^{\left[nt\right]} \mathbf{J}_i, n^{-1/\alpha}\sum_{i=1}^{\left[nt\right]}T_i \right) \overunderset{\mathbb{J}_1}{n \to \infty}{\to } \left(\mathbf{L}_{\alpha}\left(t\right), S_{\alpha}\left(t\right)\right),
        \end{equation}
        where $\{\left(\mathbf{L}_{\alpha}\left(t\right), S_{\alpha}\left(t\right)\right)\}_{t\ge0}$ is the $\left(d+1\right)$-dimensional L\'evy process described by the L\'evy triplet
        \begin{equation}\label{levymeasure1}
            \left[\int_{\mathbf{x}\neq\mathbf{0}}\frac{\mathbf{x}}{1+\left\|\mathbf{x}\right\|^2}\nu_{\left(\mathbf{L}_{\alpha}, S_{\alpha}\right)}\left(d\mathbf{x}\right), \mathbf{0}, \nu_{\left(\mathbf{L}_{\alpha}, S_{\alpha}\right)}\left(d\mathbf{x}\right)= \int_{\mathbf{D}} \int_{\mathbb{R}_+} \delta_{vt\mathbf{u}}(d\mathbf{x}) \, P_V(dv) \, \Lambda(d\mathbf{u}) \, \nu_\alpha(dt)\right],
        \end{equation}
    \item 
    for $\beta<\alpha$, we have
        \begin{equation}\label{convcase2}
            \left(n^{-1/\beta}\sum_{i=1}^{\left[nt\right]} \mathbf{J}_i, n^{-1/\alpha}\sum_{i=1}^{\left[nt\right]}T_i \right) \overunderset{\mathbb{J}_1}{n \to \infty}{\to } \left(\mathbf{L}_{\beta}\left(t\right), S_{\alpha}\left(t\right)\right),
        \end{equation}
        where $\{\left(\mathbf{L}_{\beta}\left(t\right), S_{\alpha}\left(t\right)\right)\}_{t\geq0}$ is the $\left(d+1\right)$-dimensional L\'evy process described by the L\'evy triplet
        \begin{equation}\label{levymeasure2}
        \begin{aligned}
            \left[\int_{\mathbf{x}\neq\mathbf{0}}\frac{\mathbf{x}}{1+\left\|\mathbf{x}\right\|^2}\nu_{\left(\mathbf{L}_{\beta}, S_{\alpha}\right)}\left(d\mathbf{x}\right), \mathbf{0}, \nu_{\left(\mathbf{L}_{\beta}, S_{\alpha}\right)}\left(d\mathbf{x}\right)=\right. \\ \left.
            \int_{\mathbf{D}} \left(\delta_{0}\left(dt\right)\nu_{L_{\beta}}\left(d\left(\mathbf{u}t\right)\right) + \delta_{0}\left(d\left(\mathbf{u}t\right)\right) \nu_\alpha\left(dt\right)\right)\mathbf{\Lambda}\left(d\mathbf{u}\right)\right],
        \end{aligned}
        \end{equation}
        \item 
        for $\alpha=\beta$, we have
        \begin{equation}\label{convcase3}
            \left((n\log n)^{-1/\alpha}\sum_{i=1}^{\left[nt\right]} \mathbf{J}_i, n^{-1/\alpha}\sum_{i=1}^{\left[nt\right]}T_i \right) \overunderset{\mathbb{J}_1}{n \to \infty}{\to } \left(\mathbf{L}_{\alpha}^{\text{log}}\left(t\right), S_{\alpha}\left(t\right)\right),
        \end{equation}
        where $\{\left(\mathbf{L}_{\alpha}^{\text{log}}\left(t\right), S_{\alpha}\left(t\right)\right)\}_{t\geq0}$ is the $\left(d+1\right)$-dimensional L\'evy process described by the L\'evy triplet
        \begin{equation}
        \begin{aligned}\label{levymeasure3}
            \left[\int_{\mathbf{x}\neq\mathbf{0}}\frac{\mathbf{x}}{1+\left\|\mathbf{x}\right\|^2}\nu_{\left(\mathbf{L}_{\alpha}, S_{\alpha}\right)}\left(d\mathbf{x}\right), \mathbf{0}, \nu_{\left(\mathbf{L}_{\alpha}^{\text{log}}, S_{\alpha}\right)}\left(d\mathbf{x}\right)=\right. \\ \left.\int_{\mathbf{D}} \left(\delta_{0}\left(dt\right)\nu_{L_{\alpha}}\left(d\left(\mathbf{u}t\right)\right) + \delta_{0}\left(d\left(\mathbf{u}t\right)\right) \nu_\alpha\left(dt\right)\right)\mathbf{\Lambda}\left(d\mathbf{u}\right)\right],
        \end{aligned}
        \end{equation}
    \end{enumerate}
        and $v\in \mathbb{R}_+$, $P_V\left(dv\right) = \mathbb{P}\left(V_1\in dv\right)$, $\mathbf{x} = \left(\mathbf{x}_1, t\right)\in \mathbb{R}^d\setminus\left\{\mathbf{0}\right\}\times \mathbb{R}_+$, $\mathbf{u} = \mathbf{x}_1/\left\|\mathbf{x}_1\right\|\in \mathbf{D}\subset\mathbb{S}^{d-1}$. Here $\mathbf{\Lambda}\left(d\mathbf{u}\right) = \mathbb{P}\left(\mathbf{I}_1\in d\mathbf{u}\right)$ is called the spectral measure.
    \begin{proof}
        \begin{enumerate}
            \item 
        For any Borel sets $\mathbf{B}_1 \in \mathcal{B}(\mathbb{R}^d)$ and $\mathbf{B}_2 \in \mathcal{B}(\mathbb{R}_+)$ such that $\mathbf{B}_1(R, \mathbf{D}) = \{r\mathbf{u} \in \mathbb{R}^d : r \in R, \mathbf{u} \in \mathbf{D}\}$ with $R \in \mathcal{B}(\mathbb{R}_+)$ and $\mathbf{D} \in \mathcal{B}(\mathbb{S}^{d-1})$, we have
\begin{equation}
\begin{aligned}
    &n\mathbb{P}\left(n^{-1/\alpha}V_1 T_1 \mathbf{I}_1 \in \mathbf{B}_1, \, n^{-1/\alpha}T_1 \in \mathbf{B}_2\right) \\
    &\quad = n \int_{\mathbf{D}} \mathbb{P}\left(n^{-1/\alpha}V_1 T_1 \mathbf{u} \in \mathbf{B}_1, \, n^{-1/\alpha}T_1 \in \mathbf{B}_2\right) \Lambda(d\mathbf{u}) \\
    &\quad = n \int_{\mathbf{D}} \int_{n^{1/\alpha}\mathbf{B}_2} \int_{\mathbb{R}_+} \mathbf{1}\left(n^{-1/\alpha} v t \mathbf{u} \in \mathbf{B}_1\right) P_V(dv) \, \mathbb{P}(T_1 \in dt) \, \Lambda(d\mathbf{u}) \\
    &\quad = \int_{\mathbf{D}} \int_{\mathbb{R}_+} \int_{\mathbb{R}_+} \mathbf{1}\left(v s \mathbf{u} \in \mathbf{B}_1\right) \mathbf{1}\left(s \in \mathbf{B}_2\right) P_V(dv) \, n\mathbb{P}(n^{-1/\alpha}T_1 \in ds) \, \Lambda(d\mathbf{u}).
\end{aligned}
\end{equation}
Thus, for $(\mathbf{x}, t) \in (\mathbb{R}^d \setminus \{\mathbf{0}\}) \times \mathbb{R}_+$ and $\mathbf{u} \in \mathbb{S}^{d-1}$, the measure can be written in terms of the Dirac delta
\begin{equation}
    n\mathbb{P}\left(n^{-1/\alpha}\mathbf{J}_1 \in d\mathbf{x}, \, n^{-1/\alpha}T_1 \in dt\right) = \int_{\mathbf{D}} \int_{\mathbb{R}_+} \delta_{vt\mathbf{u}}(d\mathbf{x}) \, P_V(dv) \, \Lambda(d\mathbf{u}) \, n\mathbb{P}(n^{-1/\alpha}T_1 \in dt).
\end{equation}
By letting $n \to \infty$ and using the regular variation of $T$, we obtain
\begin{equation}
    n\mathbb{P}\left(n^{-1/\alpha}\mathbf{J}_1 \in d\mathbf{x}, \, n^{-1/\alpha}T_1 \in dt\right) \xrightarrow{v} \int_{\mathbf{D}} \int_{\mathbb{R}_+} \delta_{vt\mathbf{u}}(d\mathbf{x}) \, P_V(dv) \, \Lambda(d\mathbf{u}) \, \nu_\alpha(dt).
\end{equation}
Finally, the limiting joint Lévy measure $\nu_{(\mathbf{L}_{\alpha}, S_{\alpha})}$ is given by
\begin{equation}\label{lm1}
    \nu_{(\mathbf{L}_{\alpha}, S_{\alpha})}(d\mathbf{x}, dt) = \int_{\mathbf{D}} \int_{\mathbb{R}_+} \delta_{vt\mathbf{u}}(d\mathbf{x}) \, P_V(dv) \, \Lambda(d\mathbf{u}) \, \nu_\alpha(dt).
\end{equation}
        Therefore, condition (a) of Theorem 3.2.2 \cite{meerschaert2001limit} is fulfilled. Moreover, since the L\'evy measure in \eqref{lm1} is a L\'evy measure of $d+1$-dimensional multidimensional $\alpha$-stable variable then the condition (b) of Theorem 3.2.2 \cite{meerschaert2001limit} is fulfilled with $Q_{(\mathbf{L}_\alpha,S_\alpha)}=  0$, which completes the proof of convergence in distribution. The last property together with Theorem 4.1 \cite{meerschaert2001limit} (invariance principle) implies the functional convergence \eqref{convcase1}.
        \item
            From Lemma \ref{auxlemma1} and Lemma \ref{vectorconv} we see that components $\mathbf{L}_{\beta}\left(t\right)$ and $S_{\alpha}\left(t\right)$ are independent, since sequences $\left\{V_i\right\}_{i\geq 1}$ and $\left\{T_i\right\}_{i\geq 1}$ are independent. Thus, from Corollary 2.3 in \cite{BeckerKern2004} we immediately read the form of the L\'evy triplet of the limiting process. 
        \item 
        Consider Borel sets $\mathbf{B}_1 \in \mathcal{B}(\mathbb{R}^d \setminus \{\mathbf{0}\})$ and $\mathbf{B}_2 \in \mathcal{B}((0, \infty))$. Let $\mathbf{B}_1(R, \mathbf{D}) = \{r\mathbf{u} \in \mathbb{R}^d : r \in R, \mathbf{u} \in \mathbf{D}\}$, where $R \in \mathcal{B}(\mathbb{R}_+)$ and $\mathbf{D} \in \mathcal{B}(\mathbb{S}^{d-1})$. Let us denote $a_n=(n\log n)^{1/\alpha}$, $b_n=n^{1/\alpha}$. The condition $a_n^{-1}\mathbf{J}_1 \in \mathbf{B}_1$ requires its radial component (its norm) to belong to $R$. Since $\mathbf{I}_1$ is a unit vector ($\|\mathbf{I}_1\| = 1$) and $V_1, T_1$ are positive, the norm is $\|\mathbf{J}_1\| = V_1 T_1$.  Therefore, $\|a_n^{-1}\mathbf{J}_1\| \in R \implies a_n^{-1}V_1 T_1 > x \implies V_1 T_1 > x a_n$. Similarly, $b_n^{-1}T_1 \in \mathbf{B}_2 \implies T_1 > y b_n$. Thus, we can bound the measure of the joint set by the joint tail probability 
\begin{align}
    n\mathbb{P}\left(a_n^{-1}\mathbf{J}_1 \in \mathbf{B}_1, \, b_n^{-1}T_1 \in \mathbf{B}_2\right) &\le n\mathbb{P}\left(V_1 T_1 > x a_n, \, T_1 > y b_n\right) := P_n(x, y).
\end{align}
Using the law of total probability and the independence of $V_1$ and $T_1$, we condition on $T_1 = s$. We, have
\begin{equation}
    P_n(x, y) = n \int_{y b_n}^\infty \mathbb{P}\left(V_1 > \frac{x a_n}{s}\right) \, dF_T(s) = n \int_{y b_n}^\infty \bar{F}_V\left(\frac{x a_n}{s}\right) \, dF_T(s).
\end{equation}
For any arbitrarily small $\epsilon > 0$, there exists a large threshold $M > 0$ such that for all $v \ge M$, $\bar{F}_V(v) \le (1+\epsilon)v^{-\alpha}$. We split the integral into two according to $M$
\begin{equation}
    P_n(x,y) = n \int_{y b_n}^{x a_n / M} \bar{F}_V\left(\frac{x a_n}{s}\right) \, dF_T(s) + n \int_{x a_n / M}^\infty \bar{F}_V\left(\frac{x a_n}{s}\right) \, dF_T(s) = I_1(n) + I_2(n).
\end{equation}
Using the trivial probability bound $\bar{F}_V(\cdot) \le 1$, we have
\begin{equation}
    I_2(n) \le n \int_{x a_n / M}^\infty 1 \, dF_T(s) = n \bar{F}_T\left(\frac{x a_n}{M}\right).
\end{equation}
Since $\bar{F}_T(t) \sim t^{-\alpha}$ and $a_n = (n \log n)^{1/\alpha}$, we obtain
\begin{equation}
    I_2(n) \lesssim n \left( \frac{x a_n}{M} \right)^{-\alpha} = n M^\alpha x^{-\alpha} (n \log n)^{-1} = M^\alpha x^{-\alpha} \frac{1}{\log n}.
\end{equation}
As $n \to \infty$, $I_2(n) \to 0$. For $I_1(n)$, the domain of integration guarantees the argument $v = \frac{x a_n}{s} \ge M$. Applying upper bound $\bar{F}_V(v) \le (1+\epsilon)v^{-\alpha}$, we have
\begin{equation}
    I_1(n) \le n \int_{y b_n}^{x a_n / M} (1+\epsilon) \left( \frac{x a_n}{s} \right)^{-\alpha} \, dF_T(s) = n (1+\epsilon) (x a_n)^{-\alpha} \int_{y b_n}^{x a_n / M} s^\alpha \, dF_T(s).
\end{equation}
Let $A = y b_n$ and $B = x a_n / M$. Using integration by parts and noting that $dF_T(s) = -d\bar{F}_T(s)$, we evaluate the integral
\begin{equation}
\begin{aligned}
    J(A, B) &= -\int_A^B s^\alpha \, d\bar{F}_T(s) = -\left( \left[ s^\alpha \bar{F}_T(s) \right]_A^B - \int_A^B \bar{F}_T(s) \, d(s^\alpha) \right) \\
            &= A^\alpha \bar{F}_T(A) - B^\alpha \bar{F}_T(B) + \int_A^B \bar{F}_T(s) \alpha s^{\alpha-1} \, ds.
\end{aligned}
\end{equation}
As $n \to \infty$, the boundary terms $A^\alpha \bar{F}_T(A)$ and $B^\alpha \bar{F}_T(B)$ converge to $1$ (since $\bar{F}_T(s) \sim s^{-\alpha}$), meaning their difference is bounded by a constant. The integral term dominates. Because $\bar{F}_T(s) \sim s^{-\alpha}$, the integrand is asymptotically $\alpha s^{-1}$. Thus, we have
\begin{equation}
    \int_A^B \alpha s^{-1} \, ds = \alpha (\log(B) - \log(A)) = \alpha \log\left( \frac{B}{A} \right).
\end{equation}
Substituting $A = y n^{1/\alpha}$ and $B = \frac{x}{M} (n \log n)^{1/\alpha}$, we evaluate the logarithmic term
\begin{equation}
    \alpha \log\left( \frac{\frac{x}{M} (n \log n)^{1/\alpha}}{y n^{1/\alpha}} \right) = \alpha \log\left( \frac{x}{M y} (\log n)^{1/\alpha} \right) = \alpha \log\left( \frac{x}{M y} \right) + \log(\log n) \sim \log(\log n).
\end{equation}
Therefore, $J(A, B) \sim \log(\log n)$. Substituting this asymptotic behavior back into our bound for $I_1(n)$ and using $n(a_n)^{-\alpha} = (\log n)^{-1}$, we have
\begin{equation}
    I_1(n) \lesssim (1+\epsilon) x^{-\alpha} \frac{1}{\log n} \log(\log n) = (1+\epsilon) x^{-\alpha} \frac{\log(\log n)}{\log n}.
\end{equation}
As $n \to \infty$, $I_1(n) \to 0$.
Since both $I_1(n) \to 0$ and $I_2(n) \to 0$, thus $P_n(x, y) \to 0$, components of limiting vector are independent. Thus, from Corollary 2.3 in \cite{BeckerKern2004} we immediately read the form of the L\'evy triplet of the limiting process. 
        \end{enumerate}
    \end{proof}
\end{lemma}
Now, let us state the main results of this section, which are scaling limits of L\'evy walks with random velocities. Before we state such results, we need to introduce additional notation. First, we define left-continuous version of right-continuous process $\{X\left(t\right)\}_{t\ge0}$ by putting $X^-\left(t\right) = \lim_{s\to t^-}X\left(s\right)$ and  the right-continuous version of any left-continuous process $\{Y\left(t\right)\}_t\geq0$ by putting $Y^+\left(t\right) = \lim_{s\to t^+}Y\left(s\right)$. Now, we prove a theorem describing scaling limits of LWs with random velocities.
\begin{theorem}
    Let $\{U_V\left(t\right)\}_{t\ge0}$ and $\{O_V\left(t\right)\}_{t\ge0}$ be the wait-first LW with
 random velocities and jump-first LW with random velocities, respectively, that are generated by the sequence $\left\{V_i, \, T_i, \, \mathbf{J}_i\right\}_{i\in \mathbb{N}}$. Then,
 \begin{enumerate}
     \item 
     if $\alpha<\beta$, we have
     \begin{equation}
         \begin{aligned}
             n^{-1/\alpha}U_V\left(n^{1/\alpha}t\right) \overunderset{\mathbb{J}_1}{n \to \infty}{\to } \left(\mathbf{L}_{\alpha}^-\left(S_{\alpha}^{-1}\left(t\right)\right)\right)^+,\\
             n^{-1/\alpha}O_V\left(n^{1/\alpha}t\right) \overunderset{\mathbb{J}_1}{n \to \infty}{\to } \mathbf{L}_{\alpha}\left(S_{\alpha}^{-1}\left(t\right)\right),
         \end{aligned}
     \end{equation}
     where the processes $\{\mathbf{L}_{\alpha}\left(t\right)\}_{t\ge0}$ and $\{S_{\alpha}\left(t\right)\}_{t\ge0}$ are dependent according to their joint L\'evy measure given by Eq. \eqref{levymeasure1},
     \item 
     if $\alpha>\beta$, we have
     \begin{equation}
         \begin{aligned}
             n^{-1/\beta}U_V\left(n^{1/\alpha}t\right) \overunderset{\mathbb{J}_1}{n \to \infty}{\to } \left(\mathbf{L}_{\beta}^-\left(S_{\alpha}^{-1}\left(t\right)\right)\right)^+,\\
             n^{-1/\beta}O_V\left(n^{1/\alpha}t\right) \overunderset{\mathbb{J}_1}{n \to \infty}{\to } \mathbf{L}_{\beta}\left(S_{\alpha}^{-1}\left(t\right)\right),
         \end{aligned}
     \end{equation}
     where the processes $\{\mathbf{L}_{\beta}\left(t\right)\}_{t\ge0}$ and $\{S_{\alpha}\left(t\right)\}_{t\ge0}$ are independent according their joint L\'evy measure given by Eq. \eqref{levymeasure2},
     \item 
     if $\alpha=\beta$, we have
     \begin{equation}
         \begin{aligned}
             n^{-1/\alpha}U_V\left(n^{1/\alpha}t\right) \overunderset{\mathbb{J}_1}{n \to \infty}{\to } \left(\mathbf{L}_{\alpha}^{\text{log}-}\left(S_{\alpha}^{-1}\left(t\right)\right)\right)^+,\\
             n^{-1/\alpha}O_V\left(n^{1/\alpha}t\right) \overunderset{\mathbb{J}_1}{n \to \infty}{\to } \mathbf{L}_{\alpha}^{\text{log}}\left(S_{\alpha}^{-1}\left(t\right)\right),
         \end{aligned}
     \end{equation}
     where the processes $\{\mathbf{L}_{\alpha}^{\text{log}}\left(t\right)\}_{t\ge0}$ and $\{S_{\alpha}\left(t\right)\}_{t\ge0}$ are independent according their joint L\'evy measure given by Eq. \eqref{levymeasure3}.
 \end{enumerate}
 \begin{proof}
     To prove part 1 and 3, we consider the array of random vectors $\left\{n^{-1/\alpha}\mathbf{J}_i, \, n^{-1/\alpha}T_i\right\}_{i\in\mathbb{N}}$ and define two auxiliary sequences of wait-first LW and jump-first LW
     \begin{equation}
         U_{V,\,n}\left(t\right) = n^{-1/\alpha}\sum_{i=1}^{N_n\left(t\right)}\mathbf{J}_i, \quad O_{V,\,n}\left(t\right) = n^{-1/\alpha}\sum_{i=1}^{N_n\left(t\right)+1}\mathbf{J}_i,
     \end{equation}
     where 
     \begin{equation}
         N_n\left(t\right) = \max\left\{k\geq0:\, \sum_{i=1}^k n^{-1/\alpha}T_i\leq t\right\}.
     \end{equation}
    Since Lemma \ref{vectorconv} holds, by Theorem 3.6 from \cite{Straka2011}, the process $\{S_{\alpha}\left(t\right)\}_{t\ge0}$ has strictly increasing realizations for parts 1 and 3 holds
    \begin{equation}
        \begin{aligned}
             U_{V, n}\left(n^{1/\alpha}t\right) \overunderset{\mathbb{J}_1}{n \to \infty}{\to } \left(\mathbf{L}_{\alpha}^-\left(S_{\alpha}^{-1}\left(t\right)\right)\right)^+,\quad 
             O_{V, n}\left(n^{1/\alpha}t\right) \overunderset{\mathbb{J}_1}{n \to \infty}{\to } \mathbf{L}_{\alpha}\left(S_{\alpha}^{-1}\left(t\right)\right).
         \end{aligned}
    \end{equation}
    Finally, we observe, that 
    \begin{equation}
        N_n\left(t\right) = \max\left\{k\geq0:\, \sum_{i=1}^k n^{-1/\alpha}T_i\leq t\right\} = \max\left\{k\geq0:\, \sum_{i=1}^k T_i\leq n^{1/\alpha}t\right\} = N\left(n^{1/\alpha}t\right).
    \end{equation}
    Hence, $U_{V,\,n}\left(t\right) = n^{-1/\alpha}U_V\left(n^{-1/\alpha}t\right)$ and $O_{V,\,n}\left(t\right) = n^{-1/\alpha}O_V\left(n^{-1/\alpha}t\right)$, which completes the proof of parts 1 and 3. The proof of part 2 is analogous. The L\'evy triplets of the limiting processes follows directly from Lemma \ref{vectorconv}.
 \end{proof}
\end{theorem}
\begin{theorem}
    Let $\{W_V\left(t\right)\}_{t\ge0}$ be the continuous LW with random velocities generated by the sequence $\left\{V_i, \, T_i, \, \mathbf{J}_i\right\}_{i\in \mathbb{N}}$. Then,
\begin{enumerate}
    \item if $\alpha<\beta$, we have
    \begin{equation}
        \frac{1}{n}W(nt) \xrightarrow[n \to \infty]{U} L_{V, \alpha}(t) = 
        \begin{cases}
            \big(\mathbf{L}_{\alpha}^-(S_{\alpha}^{-1}(t))\big)^+, & t\in\mathcal{R}(S_{\alpha}),\\[1ex]
            \begin{aligned}
                &\big(\mathbf{L}_{\alpha}^-(S_{\alpha}^{-1}(t))\big)^+ + \frac{t-G(t)}{H(t) - G(t)}\\
                &\quad \times \big(\mathbf{L}_{\alpha}(S_{\alpha}^{-1}(t)) - \big(\mathbf{L}_{\alpha}^-(S_{\alpha}^{-1}(t))\big)^+\big),
            \end{aligned} & t\notin \mathcal{R}(S_{\alpha}).
        \end{cases}
    \end{equation}
    
    \item if $\alpha>\beta$, we have
    \begin{equation}
        \frac{1}{n}W(nt) \xrightarrow[n \to \infty]{U} L_{V, \beta}(t) = 
        \begin{cases}
            \big(\mathbf{L}_{\beta}^-(S_{\alpha}^{-1}(t))\big)^+, & t\in\mathcal{R}(S_{\alpha}),\\[1ex]
            \begin{aligned}
                &\big(\mathbf{L}_{\beta}^-(S_{\alpha}^{-1}(t))\big)^+ + \frac{t-G(t)}{H(t) - G(t)}\\
                &\quad \times \big(\mathbf{L}_{\beta}(S_{\alpha}^{-1}(t)) - \big(\mathbf{L}_{\beta}^-(S_{\alpha}^{-1}(t))\big)^+\big),
            \end{aligned} & t\notin \mathcal{R}(S_{\alpha}).
        \end{cases}
    \end{equation}
    
    \item if $\alpha=\beta$, we have
    \begin{equation}
        \frac{1}{n}W(nt) \xrightarrow[n \to \infty]{U} L_{V, \alpha}^{\text{log}}(t) = 
        \begin{cases}
            \big(\mathbf{L}_{\alpha}^{\text{log}-}(S_{\alpha}^{-1}(t))\big)^+, & t\in\mathcal{R}(S_{\alpha}),\\[1ex]
            \begin{aligned}
                &\big(\mathbf{L}_{\alpha}^{\text{log}-}(S_{\alpha}^{-1}(t))\big)^+ + \frac{t-G(t)}{H(t) - G(t)}\\
                &\quad \times \big(\mathbf{L}_{\alpha}^{\text{log}}(S_{\alpha}^{-1}(t)) - \big(\mathbf{L}_{\alpha}^{\text{log}-}(S_{\alpha}^{-1}(t))\big)^+\big),
            \end{aligned} & t\notin \mathcal{R}(S_{\alpha}),
        \end{cases}
    \end{equation}
    where the processes 
        where $\mathcal{R}\left(S_{\alpha}\right) = \left\{S_{\alpha}\left(t\right):\, t\geq0\right\} \subset\left[0, \infty\right)$, $G\left(t\right) = S_{\alpha}^{-}\left(S_{\alpha}^{-1+}\left(t\right)\right)$, $H\left(t\right) = S_{\alpha}\left(S_{\alpha}^{-1+}\left(t\right)\right)$ and $\overunderset{U}{n \to \infty}{\to }$ means uniform convergence.
\end{enumerate}
    \begin{proof}
        According to Lemma \ref{vectorconv}, we have the convergence of the underlying joint vector of jumps and waiting times. Then, by Theorem 4.11 in \cite{Magdziarz2015b} for part 1 and 3, we have that the process
        \begin{equation}\label{scaledw}
            n^{-1/\alpha}W\left(n^{1/\alpha}t\right)
        \end{equation}
    converges weakly to $\left\{L_{V, \alpha}\left(t\right)\right\}_{t\geq0}$ in Skorokhod $\mathbb{M}_1$ topology. Because the sample paths of the limiting process $\left\{L_{V, \alpha}\left(t\right)\right\}_{t\geq0}$ are continuous according to Lemma 4.4 in \cite{Magdziarz2015b}, thus convergence also holds in the stronger uniform topology $U$. Since $b^{-1}\left(n\right) = n^{-1/\alpha}$ is regularly varying with index $1/\alpha$ there exists a function $\overset{\sim}{b}\left(n\right)$ regularly varying with index $\alpha$ such that $b\left(\overset{\sim}{b}\left(n\right)\right)^{-1} \sim n$ as $n\to\infty$. Replacing $n$ by $\overset{\sim}{b}\left(n\right)$ in Eq. \eqref{scaledw} we obtain the thesis. The proofs of the parts 2 and 3 are analogous. The form of the L\'evy triplets of the limiting processes follows directly from Lemma \ref{vectorconv}.
    \end{proof}
\end{theorem}

\section*{Acknowledgment}
M.A.T. and H.W. were supported by the National Science Centre, Poland (NCN) under the grant Sonata Bis with a number NCN 2020/38/E/ST1/00153.

\bibliographystyle{plain}
\bibliography{bibliography}
\end{document}